\documentclass{amsart}

\usepackage{amscd,amsmath,amssymb,amsfonts,verbatim}
\usepackage{times}
\usepackage[cmtip, all]{xy}
\usepackage{graphicx}
\usepackage[active]{srcltx}
\usepackage{color}
\usepackage{textcomp}

\definecolor{refkey}{gray}{.5}   
\definecolor{labelkey}{gray}{.5} 
\definecolor{Red}{rgb}{1,0,0}

\newcommand{\pf}{{\bf Proof : }}

\newcommand{\qedwhite}{\hfill \ensuremath{\Box}}

\newtheorem{theo}{Theorem}[section]
\newtheorem{prop}[theo]{Proposition}
\newtheorem{lem}[theo]{Lemma}
\newtheorem{cor}[theo]{Corollary}

\newtheorem{obs}[theo]{Observation}

\theoremstyle{definition}

\newtheorem{notn}[theo]{Notation}
\newtheorem{rem}[theo]{Remark}

\newtheorem{defi}[theo]{Definition}

\newcommand{\Um}{\mbox{\rm Um}}		\newcommand{\SL}{\mbox{\rm SL}}
\newcommand{\GL}{\mbox{\rm GL}}		
	\newcommand{\E}{\mbox{\rm E}}
\newcommand{\ESp}{\mbox{\rm ESp}}     \newcommand{\Sp}{\mbox{\rm Sp}}
      \newcommand{\EO}{\mbox{\rm EO}}

\title{ Generalised Homotopy and Commutativity Principle}
\author{Ravi A. Rao and Sampat Sharma}
\newcommand{\Addresses}{{
  \bigskip
  \footnotesize
  
  \textsc{Ravi A. Rao, School of Mathematical Sciences, Narsi Monji Institute of Management Sciences (NMIMS), 
           V.L. Mehta Road, Vile Parle West, Mumbai 400056, India}\par\nopagebreak
  \textit{E-mail:} Ravi A. Rao \texttt{<ravirao.oarivar@gmail.com>}

  \medskip
  
\textsc{Sampat Sharma, School of Mathematical and Statistical Sciences, IIT Mandi, Mandi 175005 (H.P), India}\par\nopagebreak
  \textit{E-mail:} Sampat Sharma \texttt{<sampat@iitmandi.ac.in; sampat.iiserm@gmail.com>}

  \medskip

  }}

\begin{document}

\maketitle
\subjclass 2020 Mathematics Subject Classification:{13C10, 13H99, 19B14.}

 \keywords {Keywords:}~ {Homotopy, Classical ~groups, Local-Global Principle}
 \begin{abstract}
  In this paper, we study the action of special $n\times n $ linear (resp. symplectic) matrices which are homotopic to identity on the right invertible $n\times m$ matrices. We also prove that the commutator subgroup of $\rm{O}_{2n}(R[X])$ is two stably elementary orthogonal for a local ring $R$ with $\frac{1}{2}\in R$ and $n\geq 3.$
 \end{abstract}
 
 \vskip0.50in

\begin{flushleft}
 Throughout this article we will assume $R$ to be a commutative ring with $1 \neq 0 .$
\end{flushleft}

\section{Introduction}

In $(${\cite[Corollary 1.4]{4}}$)$, Suslin established the normality of the elementary 
linear subgroup $\E_n(R)$ in $\GL_n(R)$, for $n \geq 3$. This was a major 
surprise at that time as it was known due to the work of Cohn in 
\cite{cohn} that in general $\E_2(R)$ is not normal in $\GL_2(R)$. This is 
the initial precursor to study the non-stable $K_1$ groups ${\SL_n(R)}/{\E_n(R)},~n\geq 3$.

This theorem can also be got as a consequence of the local-global principle of 
 Quillen (for projective modules) in \cite{qu}; and its analogue for the 
linear group of elementary matrices $\E_n(R[X])$, when $n \geq 3$ due to Suslin in \cite{4}. In fact, in \cite{22} it is shown that, in some sense,
the normality property of the elementary group $\E_n(R)$ in $SL_n(R)$ is 
equivalent to having a local-global principle for $\E_n(R[X])$.  

In {\cite{bak}}, Bak proved the following beautiful result:

\begin{theo} $($Bak$)$ For an almost commutative ring $R$ with identity with 
centre $C(R)$. The group ${\SL_{n}(R)}/{\E_{n}(R)}$ 
 is nilpotent of class atmost $\delta(C(R)) + 3 - n$, where  
$\delta(C(R)) < \infty ~and~ n\geq 3$, where  $\delta(C(R))$ is 
the Bass--Serre- dimension of $C(R)$.
\end{theo}

This theorem, which is proved by a localization and completion technique, 
which evolved from an adaptation of the proof of the 
Suslin's $K_1$-analogue of Quillen's  
local-global principle was further investigated in \cite{rs}. In \cite{rs}, we proved that

\begin{theo} 
\label{ab}
Let $R$ be a local ring, and let $A = R[X]$. Then the 
group ${\SL_{n}(A)}/{\E_{n}(A)}$ is an abelian group for $n \geq 3$. 
\end{theo} 

This theorem is a simple consequence of the following principle:

\vskip0.15in

\begin{theo} 
\label{homcom} $(${\cite[Theorem 2.19]{rs}}$)$ $($Homotopy and commutativity principle$):$ 
Let $R$ be a commutative ring. Let $\alpha \in \SL_n(R)$, $n \geq 
3$, be homotopic to the identity. Then, for any $\beta \in \SL_n(R)$, 
$\alpha \beta = \beta \alpha \varepsilon$, for some $\varepsilon \in \E_n(R)$. 
\end{theo}

\vskip0.15in

This principle is a consequence of the Quillen--Suslin's local-global principle; 
and using a {\it non-symmetric} application of it as done by Bak in 
\cite{bak}. 

\vskip0.15in
Using Bak's localization method, in \cite{stepnov}, Stepanov proved the following result for all simply connected Chevalley group of rank $> 1$  :
\begin{theo} Let $G$ be a simply connected group of rank $> 1$ with $G(R) = E(R)$ when $R$ is a local ring. Then for any commutative ring $R$ with 1, 
$$ [\overset{\sim}{E}(R), G(R)] = E(R)$$ 
where $\overset{\sim}{E}(R) = \bigcap_{(s_{1}, \ldots s_{l})\in Um_{l}(R)}\prod G(R, s_{1}R), \ldots, G(R, s_{l}R),$ denotes the extended elementary group. 
\end{theo}

\vskip0.15in
 In this paper, we generalise the homotopy and commutativity principle to any $n\times m$ right invertible matrix over a commutative ring 
$R$. In particular, we prove that : 

\begin{theo}$($Generalised homotopy and commutativity principle$)$
Let $R$ be a commutative ring and $V\in \Um_{n,m}(R)$ with $m> n \geq 2 ~ \mbox{or}~ m = n\geq3.$ 
 Let $\delta \in \SL_{n}(R)$ be homotopic to identity. 
 Let $\delta(T)$ be a homotopy of $\delta.$
 Then $\exists ~\sigma(T) \in \SL_{m}(R[T],(T))$ such that
 $$\delta(T) V = V\sigma(T)~ \mbox{and}~ \sigma(T)^{-1}(\delta(T) \perp I_{m-n}) \in \E_{m}(R[T],(T)).$$
 Moreover, if $\sigma(1) = \sigma$, then we have
 $\delta V = V\sigma$ and $\sigma^{-1}(\delta \perp I_{m-n}) \in \E_{m}(R).$
\end{theo}

\vskip0.15in

We also prove similar results in the case of symplectic groups (see theorem \ref{generalisedhomotopysymplectic}). We prove the similar statement in the case of orthogonal groups as well with $m\geq n+2, n\geq 2$ (see theorem \ref{generalisedhomotopyorthogonal}).
As a consequence we prove that linear and symplectic quotients are abelian, but in the 
case of orthogonal quotients we could only establish the following: 
\begin{theo}
Let $m\geq 3$, $R$ be a local ring, $\frac{1}{2}\in R.$ Then $([\rm{O}_{2m}R[X], \rm{O}_{2m}R[X]] \perp I_{2}) \subseteq \EO_{2m+2}(R[X]).$
\end{theo}

\vskip0.15in 
We do believe that orthogonal quotient groups are also abelian; as it is the case when the base ring is regular local ring 
containing a field (see $(${\cite[Corollary 4.21]{rs}}$)$.

\section{\bf{Generalised Homotopy and Commutativity Principle for Linear Groups}}

\par 
Let $v = (a_{0},a_{1},\ldots, a_{r}), w = (b_{0},b_{1},\ldots, b_{r})$ be two rows of length $r+1$ over a commutative ring 
$R$. A row $v\in R^{r+1}$ is said to be unimodular if there is a $w \in R^{r+1}$ with
$\langle v ,w\rangle = \Sigma_{i = 0}^{r} a_{i}b_{i} = 1$ and ${\Um}_{r+1}(R)$ will denote
the set of unimodular rows (over $R$) of length 
$r+1$.
\par 
The group of elementary matrices is a subgroup of $\GL_{r+1}(R)$, denoted by $\E_{r+1}(R)$, and is generated by the matrices 
of the form $e_{ij}(\lambda) = I_{r+1} + \lambda \E_{ij}$, where $\lambda \in R, ~i\neq j, ~1\leq i,j\leq r+1,~
\E_{ij} \in M_{r+1}(R)$ whose $ij^{th}$ entry is $1$ and all other entries are zero. The elementary linear group 
$\E_{r+1}(R)$ acts on the rows of length $r+1$ by right multiplication. Moreover, this action takes unimodular rows to
unimodular
rows : ${{\Um}_{r+1}(R)}/{\E_{r+1}(R)}$ will denote the set of orbits of this action; and we shall denote by $[v]$ the 
equivalence class of a row $v$ under this equivalence relation.

\begin{defi}
 An $\alpha \in M_{n \times m}(R)$ is said to be right invertible if  $\exists ~\beta \in M_{m \times n}(R)$ 
 such that $\alpha \beta = I_{n}.$ 
 We will denote set of all $ n\times m$ right invertible matrices by ${\Um}_{n,m}(R).$
\end{defi}

\begin{defi}
 An $R$-module $P$ is said to be stably free of type $n$, if $P\oplus R^{n}$ is a free module.
\end{defi}

\begin{flushleft}
 \par To every $\alpha \in {\Um}_{n,m}(R)$, we can associate a stably free module $P$ of type $n$, in the following way:\\
 \par Since $\alpha \in {\Um}_{n,m}(R)$, it gives rise to a surjective map $R^{m} \overset{\alpha} \longrightarrow R^{n}$. 
 Let $P = \mbox{Ker}({\alpha})$, then we have a short 
 exact sequence $$ 0\longrightarrow P\longrightarrow R^{m} \longrightarrow R^{n} \longrightarrow 0.$$
 Since $R^{n}$ is a free module, the above short exact sequence splits and we have $P\oplus R^{n} \simeq R^{m}.$\\
 \par To every stably free module $P$ of type $n$, we can associate an element $\alpha$ of ${\Um}_{n,m}(R)$, for some $m$, in  the
 following way:\\
 Since $P$ is stably free, we have a short exact sequence 
 $$ 0\longrightarrow P\longrightarrow R^{m} \longrightarrow R^{n} \longrightarrow 0.$$ Let $\alpha$ to be the matrix of the 
 map $R^{m}\longrightarrow R^{n}$. Since $R^{n}$ is a free module,
 the above short exact sequence splits and we have $\alpha \in {\Um}_{n,m}(R).$
\end{flushleft}

\begin{lem}
\label{stablyfreeunimodular}
 $(${\cite[Chapter 1, Proposition 4.3]{2''}}$)$ An $\alpha \in {\Um}_{n,m}(R)$ is completable to an invertible 
 matrix of determinant $1$ if and only if the corresponding stably free module is free.
\end{lem}

\begin{lem}
\label{0.4}
 Let $R$ be a local ring and $V\in {\Um}_{n,m}(R)$ for $m ~(\mbox{or}~ n)\geq 2$. Then $V$ is completable to an elementary 
 matrix.
\end{lem}
${\pf}$ Every  $V \in \mbox{Um}_{n,m}(R)$ corresponds to a stably free module $P.$ Since a projective module over a local ring is
free, $P$ is free. In view of lemma \ref{stablyfreeunimodular}, V is 
completable to a matrix $W \in \SL_{m}(R) = \E_{m}(R).$
$~~~~~~~~~~~~~~~~~~~~~~~~~~~~~~~~~~~~~~~~~~~~~~~~~~~~~~~~~~~~~~~~~~~~~~~~~~~~~~~~~~~~~~~~~~~~~~~~~~~~~~~~~~~~~~\qedwhite$

\begin{defi}
 Let R be a ring. A matrix $\alpha \in \SL_{n}(R)$ is said to be 
 homotopic to identity if there exists a matrix $\gamma(X) \in 
 \SL_{n}(R[X])$ such that $\gamma(0) = \textit{Id}~and~ \gamma(1) = \alpha$.
\end{defi}

\begin{prop}
\label{localgeneralisedhomotopy}
 Let $R$ be a local ring and $V \in \Um_{n,m}(R)$ for $m>n\geq 2$ or $m = n\geq 3.$
  Let $\delta \in \SL_{n}(R)$ be homotopic to identity. 
 Let $\delta(T)$ be a homotopy of $\delta.$ Then there exists,
 $ \sigma(T) \in \SL_{m}(R[T])$ with $\sigma(0) = \textit{Id}$ and $\sigma(T)^{-1}(\delta(T) \perp I_{m-n}) 
 \in \E_{m}(R[T])$ 
 such that $$\delta(T)V = V\sigma(T).$$
\end{prop}
${\pf}$ In view of lemma \ref{0.4}, $V$ is completable to a matrix $W \in \SL_{m}(R)$. Since $R$ is a local ring, 
$W\in \E_{m}(R).$ By $(${\cite[Corollary 1.4]{4}}$)$, $\E_{m}(R[T]) \trianglelefteq \SL_{m}(R[T])$, for $m\geq 3.$ 
Thus there exists 
$\varepsilon_{1}(T) \in \E_{m}(R[T])$ such that 
$$(\delta(T)\perp I_{m-n})W(\delta(T)\perp I_{m-n})^{-1} = \varepsilon_{1}(T).$$
 Thus we have $(\delta(T)\perp I_{m-n})W= 
\varepsilon_{1}(T)W^{-1}W(\delta(T)\perp I_{m-n}).$ Again by normality of $\E_{m}(R[T])$ in $\SL_{m}(R[T]$ for 
$m\geq 3,$ there exists $\varepsilon(T) \in \E_{m}(R[T])$  such that 
 $$(\delta(T)\perp I_{m-n})W = W(\delta(T)\perp I_{m-n})\varepsilon(T).$$
Note that $\varepsilon(0) = \textit{Id}.$
Upon taking $\sigma(T) = (\delta(T)\perp I_{m-n})\varepsilon(T)$ and multiplying above equation by$ \begin{bmatrix}
                 I_{n} & 0\\
                 0 & 0\\
                \end{bmatrix}$, we gets desired result.
$~~~~~~~~~~~~~~~~~~~~~~~~~~~~~~~~~~~~~~~~~~~~~~~~~~~~~~~~~~~~~~~~~~~~~~~~~~~~~~~~~~~~~~~~~~~~~~~~~~~~~~~~~~~~~~\qedwhite$

\begin{theo}
\label{generalisedhomotopylinear}
 Let $R$ be a commutative ring and $V\in \Um_{n,m}(R)$ with $m> n \geq 2 ~ \mbox{or}~ m = n\geq3.$ 
 Let $\delta \in \SL_{n}(R)$ be homotopic to identity. 
 Let $\delta(T)$ be a homotopy of $\delta.$
 Then $\exists ~\sigma(T) \in \SL_{m}(R[T],(T))$ such that
 $$\delta(T) V = V\sigma(T)~ \mbox{and}~ \sigma(T)^{-1}(\delta(T) \perp I_{m-n}) \in \E_{m}(R[T],(T)).$$
 Moreover
 , if $\sigma(1) = \sigma$, then we have
 $\delta V = V\sigma$ and $\sigma^{-1}(\delta \perp I_{m-n}) \in \E_{m}(R).$
\end{theo}
${\pf}$  Define,
\begin{align*}J = &\{s\in R~ | ~\delta(T)_{s}V_{s} = V_{s}\sigma(T)~ \mbox{for ~some}~
 \sigma(T) \in \SL_{m}(R_{s}[T], (T))~\\
&\mbox {with}~
\sigma(T)^{-1}(\delta(T)_{s}\perp I_{m-n}) \in \E_{m}(R_{s}[T])\}.
\end{align*}
$\bf{Claim:}$ $J$ is an ideal.
\par
 For $s\in J, \lambda \in R$, clearly $\lambda s \in J$.
 So we need to prove that if $s_{1}, s_{2} \in J$ then
 $s_{1} + s_{2} \in J.$ Since $s_{1}, s_{2} \in J$,
we have $(s_{1} + s_{2})s_{1}$, $(s_{1} + s_{2})s_{2} \in J.$
We rename $R_{s_{1}+s_{2}}$ by $R$, now
it suffices to show that
\begin{align*}\delta(T)V &= V\sigma(T) ~\mbox{for}~\mbox{some}~ \sigma(T)
\in \SL_{m}(R[T], (T))\\
~&\mbox{with}~ 
\sigma(T)^{-1}(\delta(T)\perp I_{m-n})\in \E_{m}(R[T], (T))~
\mbox{provided that}~ s_{1} + s_{2} = 1~\mbox{ and}
\end{align*}
\begin{equation}
\label{eq1}
\delta(T)_{s_{1}}V_{s_{1}} = V_{s_{1}}\sigma_{1}(T) ~\mbox{with}~
\sigma_{1}(T)^{-1}(\delta(T)_{s_{1}}\perp
I_{m-n})\in \E_{m}(R_{s_{1}}[T], (T)),
\end{equation}
\begin{equation}
\label{eq2}
 \delta(T)_{s_{2}}V_{s_{2}} = V_{s_{2}}\sigma_{2}(T) ~\mbox{with}~ 
 \sigma_{2}(T)^{-1}(\delta(T)_{s_{2}}\perp
I_{m-n})\in \E_{m}(R_{s_{2}}[T], (T)).
\end{equation}
Let  $$~~~~~~~~~~~~~~~~~~~~~~~~~~~~\sigma_{1}(T)(\delta(T)_{s_{1}} \perp
I_{m-n})^{-1} = \varepsilon_{1}(T) 
\in \E_{m}(R_{s_{1}}[T], (T)),$$
       $$~~~~~~~~~~~~~~~~~~~~~~~~~~~~\sigma_{2}(T)(\delta(T)_{s_{2}} \perp 
       I_{m-n})^{-1} = \varepsilon_{2}(T) 
       \in \E_{m}(R_{s_{2}}[T], (T)).$$

 Now, $\begin{bmatrix}
                 V_{s_{1}s_{2}}\\
                  0\\
                \end{bmatrix} \varepsilon_{1}(T)_{s_{2}}
                \varepsilon_{2}(T)_{s_{1}}^{-1}
                = \begin{bmatrix}
                 V_{s_{1}s_{2}}\\
                  0\\
                \end{bmatrix}.$ Let $\theta (T) =
                \varepsilon_{1}(T)_{s_{2}}\varepsilon_{2}(T)_{s_{1}}^{-1}.$ 
By Quillen's splitting property, $\mbox{for}~b\in (s_{2}^{N}),~ N>>0$, we have 
\begin{equation}
 \label{eq3}
 \theta(T) = \theta(bT)\{ \theta(bT)^{-1}\theta(T)\}\end{equation}
$ ~\mbox{with}
 ~\theta(bT)\in \E_{m}(R_{s_{1}}[T]), ~ \mbox{and}~ 
 \theta(bT)^{-1}\theta(T)\in  \E_{m}(R_{s_{2}}[T]).$

Since $\begin{bmatrix}
                 V_{s_{1}s_{2}}\\
                  0\\
                \end{bmatrix} \theta(T) = \begin{bmatrix}
                 V_{s_{1}s_{2}}\\
                  0\\
                \end{bmatrix}$. We have, $\begin{bmatrix}
                 V_{s_{1}s_{2}}\\
                  0\\
                \end{bmatrix} \theta(bT) = \begin{bmatrix}
                 V_{s_{1}s_{2}}\\
                  0\\
                \end{bmatrix}$ and \\ $\begin{bmatrix}
                 V_{s_{1}s_{2}}\\
                  0\\
                \end{bmatrix} \theta(bT)^{-1}\theta(T)
                 = \begin{bmatrix}
                 V_{s_{1}s_{2}}\\
                  0\\
                \end{bmatrix}.$
                \par
                Define, $\theta(bT) = \eta_{1}(T)^{-1}, ~ \theta(bT)^{-1}\theta(T) = \eta_{2}(T)$. Thus we have 
                $V_{i}\eta_{i}(T) = V_{i}$ with $\eta_{i}(0) = \textit{Id}$ and $\eta_{i}(T) \in \E_{m}(R_{s_{i}}[T])$ 
                for $i =1,2.$
                \par
                In view of equation \ref{eq3}, we gets 
                \begin{equation}
 \label{eq4}
 (\eta_{1}(T)\varepsilon_{1}(T))_{s_{2}} = (\eta_{2}(T)\varepsilon_{2}(T))_{s_{1}}.
\end{equation}
Now, by equation \ref{eq3} and equation \ref{eq4}, 
$$\delta(T)_{s_{1}}V_{s_{1}} = V_{s_{1}}\eta_{1}(T)\sigma_{1}(T)$$
$$\delta(T)_{s_{2}}V_{s_{2}} = V_{s_{2}}\eta_{2}(T)\sigma_{2}(T).$$
In view of equation 4, we have $(\eta_{1}(T)\sigma_{1}(T))_{s_{2}} = (\eta_{2}(T)\sigma_{2}(T))_{s_{1}}$. Since 
$s_{1} + s_{2} = 1$, 
$\exists ~\sigma(T) \in \SL_{m}(R[T])$ such that $\sigma(T)_{s_{1}} = \eta_{1}(T)\sigma_{1}(T)$ and 
$\sigma(T)_{s_{2}} = \eta_{2}(T)\sigma_{2}(T)$ with $\sigma(T)_{s_{i}}^{-1}(\delta(T)_{s_{i}}\perp I_{m-n}) \in \E_{m}(R_{s_{i}}
[T])$ for $i = 1, 2.$ Since $s_{1}$ and $s_{2}$ are comaximal, 
by Suslin's local-global principle  
$(${\cite[Theorem 3.1] {4}}$)$, we have
\begin{align*}\delta(T)V &= V\sigma(T) ~\mbox{for}~\mbox{some}~ 
 \sigma(T) \in \SL_{m}(R[T], (T)) ~\\
&\mbox{with}~ 
\sigma(T)^{-1}(\delta(T)\perp I_{m-n})\in \E_{m}(R[T], (T)).
\end{align*}
This proves that $J$ is an ideal.
\par
In view of Proposition \ref{localgeneralisedhomotopy}, for every maximal ideal $\mathfrak{m}$ of $R$, we have 
$$\delta(T)_{\mathfrak{m}}V_{\mathfrak{m}} = V_{\mathfrak{m}}\sigma^{'}(T) ~\mbox{with}~ 
\sigma^{'}(T)^{-1}(\delta(T)_{\mathfrak{m}}\perp
I_{m-n})\in \E_{m}(R_{\mathfrak{m}}[T], (T)).$$
Thus there exists $s \in R\setminus \mathfrak{m}$, such that
$$\delta(T)_{s}V_{s} = V_{s}\sigma^{'}(T) ~\mbox{with}~ 
\sigma^{'}(T)^{-1}(\delta(T)_{s}\perp
I_{m-n})\in \E_{m}(R_{s}[T], (T)).$$
Therefore $J \nsubseteq \mathfrak{m}$, for any maximal ideal $\mathfrak{m}$ of $R$ i.e. $1 \in J.$
Thus $\exists ~\sigma(T) \in SL_{m}(R[T], (T))$ such that 
$$\delta(T)V = V\sigma(T) ~\mbox{with}~ 
\sigma(T)^{-1}(\delta(T)\perp I_{m-n})\in \E_{m}(R[T], (T)).$$
Now put $T=1$, and take $\sigma(1) = \sigma$ to get the desired result.
$~~~~~~~~~~~~~~~~~~~~~~~~~~~~~~~~~~~~~~~~~~~~~~~~~~~~~~~~~~~~~~~~~~~~~~~~~~~~~~~~~~~~~~~~~~~~~~~~~~~~~~~~~~~~~~\qedwhite$

\begin{cor}
 $(${\cite[Theorem 2.19]{rs}}$)$
 Let $n\geq 3$ and $\alpha, \beta \in \SL_{n}(R).$ Let 
 either $\alpha$ or $\beta$ be homotopic to identity. Then $\alpha \beta = \beta \alpha \varepsilon,~
  ~\mbox{for~ some} ~\varepsilon \in \E_{n}(R). $
\end{cor}
${\pf}$ Let us assume that $\alpha$ is homotopic to identity, so there exists
 $\delta(T) \in \SL_{n}(R[T])$ such that $\delta(0) = \textit{Id}~\mbox{and}~ \delta(1) = \alpha.$ By theorem
 \ref{generalisedhomotopylinear}, 
 there exists $\varepsilon(T) \in \E_{n}(R[T])$ with $\varepsilon(0) = \textit{Id}$ such that 
 $$\delta(T)\beta = \beta\delta(T)\varepsilon(T).$$
 \par Put $T =1$ to get the desired result.
 $~~~~~~~~~~~~~~~~~~~~~~~~~~~~~~~~~~~~~~~~~~~~~~~~~~~~~~~~~~~~~~~~~~~~~~~~~~~~~~~~~~~~~~~~~~~~~~~~~~~~~~~~~~~~~~\qedwhite$

\begin{cor}
$($Vaserstein$)$
 Let $\delta \in \SL_{n}(R)$ and $V\in \Um_{n,m}(R), m> n \geq 2~\mbox{or}~n = m \geq 3.$ 
 Then $\delta V = V\sigma$ for some $\sigma \in \SL_{m}(R)$ with
  $(\sigma \perp \delta^{-1}) \in \E_{n+m}(R).$
\end{cor}

${\pf}$ By Whitehead's Lemma, $(\delta \perp \delta^{-1}) \in \E_{2m}(R).$ Since every elementary matrix is 
homotopic to identity, 
thus by theorem \ref{generalisedhomotopylinear}, 
$$(\delta \perp \delta^{-1})(V\perp I_{n}) = (V\perp I_{n})\sigma^{'}, ~\mbox{with}~\sigma^{'} \in \E_{n+m}(R).$$
Write $\sigma^{'} = \begin{bmatrix}
                 \alpha & \beta\\
                 \gamma & \zeta\\
                \end{bmatrix}$ where $\alpha \in M_{m\times m}(R),~\beta\in M_{m\times n}(R), ~ \gamma \in M_{n\times m}(R),
                ~ \zeta \in M_{n\times n}(R)$. Thus we have, 
                $$(\delta V \perp \delta^{-1}) = \begin{bmatrix}
                 V\alpha & V\beta\\
                 \gamma& \zeta\\
                \end{bmatrix}.$$
 \par
 Upon compairing both sides we gets $\gamma= 0$ and $\zeta = \delta^{-1}.$ Therefore
 $$\begin{bmatrix}
                 \alpha & \beta\\
                 0 & \delta^{-1}\\
                \end{bmatrix} \in \E_{n+m}(R).$$
\par
Now, take $\alpha = \sigma$, so we have $(\sigma \perp \delta^{-1}) \in E_{n+m}(R)$ and $\delta V = V\sigma.$
$~~~~~~~~~~~~~~~~~~~~~~~~~~~~~~~~~~~~~~~~~~~~~~~~~~~~~~~~~~~~~~~~~~~~~~~~~~~~~~~~~~~~~~~~~~~~~~~~~~~~~~~~~~~~~~\qedwhite$

\begin{lem}
\label{suslincommoperpendiclular}
 $($Suslin$)$ $(${\cite[Lemma 2.8]{sus2}}$)$ Let $r\geq 3$ and $v_{1}, ~v_{2},~ w~\in M_{1,r}(R)$ be such that 
 $\langle v_{1}, w\rangle = \langle v_{2},w\rangle = 1$, then $v_{1} \overset{\E}\sim v_{2}.$
\end{lem}

\begin{cor}
 \label{applicationsuslincommonperpendicular}
 Let $n\geq 3$ and $\alpha = \begin{bmatrix}
                 a_{1}, & \ldots, & a_{n}\\
                 b_{1}, & \ldots, & b_{n}\\
                \end{bmatrix} \in \Um_{2, n}(R).$ Then, 
                $$(a_{1}, \ldots, a_{n}) \overset{\E}\sim (b_{1}, \ldots, b_{n}).$$
\end{cor}
${\pf}$  Since $\alpha \in \mbox{Um}_{2,n}(R)$, $\exists ~ \beta = \begin{bmatrix}
                 c_{1} & d_{1}\\
                 \vdots & \vdots \\
                 c_{n} & d{_n} \\
                \end{bmatrix} \in M_{n,2}(R),$ such that $\alpha \beta = I_{2}.$ Let\\ $w = (c_{1}+d_{1}, \ldots, 
                c_{n}+d_{n})$.
                Since $\langle(a_{1}, \ldots, a_{n}),w\rangle = \langle(b_{1}, \ldots, b_{n}),w\rangle = 1.$ 
                Thus by lemma \ref{suslincommoperpendiclular}, 
                $$(a_{1}, \ldots, a_{n}) \overset{E}\sim (b_{1}, \ldots, b_{n}).$$
      $~~~~~~~~~~~~~~~~~~~~~~~~~~~~~~~~~~~~~~~~~~~~~~~~~~~~~~~~~~~~~~~~~~~~~~~~~~~~~~~~~~~~~~~~~~~~~~~~~~~~~~~~~~~~~~\qedwhite$

     \begin{cor}
     $($Roitman$)$ $(${\cite[Theorem 8]{roit}}$)$
      Let $(x_{0}, \ldots, x_{n})\in \Um_{n+1}(R),~n\geq 2$ and $0\leq k \leq n-1,~y_{i}\in R~\mbox{for}~k\leq i\leq n.$ 
      Let I be an ideal of $R$ generated by $2\times 2$ minors of the matrix
      $$ \alpha = \begin{bmatrix}
                 x_{k}, & \ldots, & x_{n}\\
                 y_{k}, & \ldots, & y_{n}\\
                \end{bmatrix}.$$
           Assume that $Rx_{0} +\ldots + Rx_{k-1} + I = R.$ Then
           $$(x_{0}, \ldots, x_{k-1}, x_{k}, \ldots, x_{n}) \overset{E}\sim (x_{0}, \ldots, x_{k-1}, y_{k},\ldots,
           y_{n})
           .$$
      \end{cor}
	${\pf}$ Consider the ring $\overset{-}R = R/Rx_{0} +\ldots + Rx_{k-1}$, by hypothesis we have
	$\overset{-}R = \overset{-}I$, therfore
	 $\overset{-}\alpha \in \mbox{Um}_{2,n-k+1}(\overset{-}R).$ Thus by corollary \ref{applicationsuslincommonperpendicular}, 
	 $\exists ~\overset{-}\varepsilon \in 
	 \E_{n-k+1}(\overset{-}R)$ such that $(\overset{-}x_{k},\ldots,\overset{-}x_{n})\overset{-}\varepsilon = 
	 (\overset{-}y_{k},\ldots,\overset{-}y_{n}).$ Let $\varepsilon \in \E_{n-k+1}(R)$ be a lift of $\overset{-}\varepsilon.$ 
	 Therefore, 
	 $$(x_{k},\ldots,x_{n})\varepsilon = (y_{k}+a_{k},\ldots,y_{n}+a_{n}),~\mbox{for~some}~a_{i}\in 
	 Rx_{0} +\ldots + Rx_{k-1}.$$
	 Thus we have, $(x_{0}, \ldots, x_{k-1}, x_{k}, \ldots, x_{n})(I_{k}\perp \varepsilon) = 
	 (x_{0}, \ldots, x_{k-1}, y_{k}+a_{k},\ldots,
           y_{n}+a_{n}).$ $\mbox{Since}~ a_{i} \in Rx_{0} +\ldots + Rx_{k-1}$, we have
           $$(x_{0}, \ldots, x_{k-1}, x_{k}, \ldots, x_{n}) \overset{E}\sim (x_{0}, \ldots, x_{k-1}, y_{k},\ldots,
           y_{n}).$$
           $~~~~~~~~~~~~~~~~~~~~~~~~~~~~~~~~~~~~~~~~~~~~~~~~~~~~~~~~~~~~~~~~~~~~~~~~~~~~~~~~~~~~~~~~~~~~~~~~~~~~~~~~~
           ~~~~~\qedwhite$

  \section{\bf{Generalised Homotopy and Commutativity Principle for Symplectic Groups}} 
  
  \begin{notn}

Let $\psi_{1} = \begin{bmatrix}
                 0 & 1\\
                 -1 & 0\\
                \end{bmatrix},~~ \psi_{n} = \psi_{n-1} \perp \psi_{1};$~~  for $n > 1$.
\end{notn}

\begin{notn}
 Let $\sigma$ be the permutation of the natural numbers given by $\sigma(2i) = 2i-1$ and $\sigma(2i-1) = 2i$.
\end{notn}

\begin{notn}
 $\E_{ij}(\lambda)$ will denote a matrix whose $ij^{th}$ entry is $\lambda$ and all other entries are $0.$
\end{notn}

\begin{defi}{\bf{Symplectic group}}~{$\Sp_{2m}(R) :$} The subgroup of $\GL_{2m}(R)$ consisting of all $2m \times 2m $ matrices 
$\{\alpha \in \GL_{2m}(R)~\mid \alpha^{t}\psi_{m}\alpha = \psi_{m}\}$.
 
\end{defi}

\begin{defi}{\bf{Elementary symplectic group}}~{$\ESp_{2m}(R)$}: We define for $ 1\leq i \neq j\leq 2m,~z\in R,$\\
$$
se_{ij}(z)=
\begin{cases}
I_{2m} + zE_{ij},~~~\textit{if}~ i = \sigma(j);\\
 I_{2m} + zE_{ij} - (-1)^{i+j}zE_{\sigma(j)\sigma(i)}, ~~~ \textit{if}~ i\neq \sigma(j).

\end{cases}
$$
\end{defi}

It is easy to verify that all these matrices belong to $\Sp_{2m}(R)$. We call them the elementary symplectic matrices
over $R$. The subgroup generated by them is called the elementary symplectic group and is denoted by $\ESp_{2m}(R)$.

\begin{notn}
 $\Sp \Um_{2n, 2m}(R) = \{V \in \Um_{2n,2m}(R) | V\psi_{m}V^{t} = \psi_{n}\}.$
\end{notn}

\begin{lem}
 $($Rao-Swan$)$ Let $n\geq 2$ and $\varepsilon \in \E_{2n}(R)$. Then there exists $\rho \in \E_{2n-1}(R)$ such that 
 $\varepsilon (1\perp \rho) \in \ESp_{2n}(R).$
\end{lem}
${\pf}$ For a proof see $(${\cite[Lemma 4.4]{pr}}$).$
$~~~~~~~~~~~~~~~~~~~~~~~~~~~~~~~~~~~~~~~~~~~~~~~~~~~~~~~~~~~~~~~~~~~~~~~~~~~~~~~~~~~~~~~~~~~~~~~~~~~~~~~~~~~~~~\qedwhite$ 

\begin{lem}
\label{vasobs}
 $($ Vaserstein$)$ $(${\cite[Lemma 5.5]{7}}$)$ For an associative ring $R$ with identity, and for any natural number $m$
 $$e_{1}\E_{2m}(R) = e_{1}(\Sp_{2m}(R) \cap \E_{2m}(R)).$$
\end{lem}

\begin{rem}
\label{vasremark}
 It was observed in $(${\cite[Lemma 2.13]{3''}}$)$  that Vaserstein's proof actually shows that  
 $e_{1}\E_{2m}(R) = e_{1}\ESp_{2m}(R).$
\end{rem}

\begin{theo}
\label{2.17}
{\bf (Local-Global principle for the symplectic groups)} 
$(${\cite[Theorem 3.6] {kop}}$)$
 Let $m\geq 2$ and $\alpha(X) \in \Sp_{2m}(R[X])$, with $\alpha(0) = \textit{Id}$. Then $\alpha(X) \in \ESp_{2m}(R[X])$ 
 if and only if for any maximal ideal $\mathfrak{m} \subset R$, the canonical image of $\alpha(X) \in 
 \Sp_{2m}(R_{\mathfrak{m}}[X])$ lies in $\ESp_{2m}(R_{\mathfrak{m}}[X])$.
\end{theo}

\begin{lem}
\label{localsymplecticcompletable}
 Let $R$ be a local ring, $m\geq n\geq 1$ and $V\in \Sp \Um_{2n,2m}(R).$ Then $V$ is completable to an elementary symplectic
 matrix.
\end{lem}

${\pf}$ We will proceed by induction on $n.$ Since $\Sp_{2}(R) = \ESp_{2}(R),$ we are done for the case $m = n =1.$
Let us assume that $n =1, m>1$,
since $V \in \Sp \Um_{2,2m}(R) \subseteq \Um_{2,2m}(R)$ and $R$ is a local ring, there exists 
$\varepsilon \in \E_{2m}(R)$ such that 
$$ V\varepsilon = \begin{bmatrix}
                 1  & 0 &  \cdots & 0\\
                 0 & 1 & \cdots & 0\\
                \end{bmatrix}.$$
\par In view of Rao-Swan Lemma, there exists $\rho \in \E_{2n-1}(R)$ such that $\varepsilon (1\perp \rho) \in \ESp_{2n}(R),$
 therefore 
 $$ V\varepsilon (1\perp \rho) = \begin{bmatrix}
                 1  & 0 &  \cdots & 0\\
                 0 & b_{2} & \cdots & b_{2m}\\
                \end{bmatrix} \in \Sp \Um_{2,2m}(R), ~\mbox{for~some}~b_{i} \in R, ~2\leq i \leq 2m.$$
Now, 
$$ \begin{bmatrix}
                 1  & 0 &  \cdots & 0\\
                 0 & b_{2} & \cdots & b_{2m}\\
                \end{bmatrix} \psi_{m} \begin{bmatrix}
                 1 & 0\\
                 0 & b_{2}\\
                 \vdots & \vdots \\
                 0 & b_{2m} \\
                \end{bmatrix} = \psi_{1}.$$
\par Upon comparing  coefficients we gets $b_{2} = 1.$ Therefore $V\overset{\ESp} \sim \begin{bmatrix}
                 1  & 0 & 0 & \cdots & 0\\
                 0 & 1 & b_{3} &\cdots & b_{2m}\\
                \end{bmatrix}.  $  Now take
                $\alpha = \prod_{k=3}^{2m}se_{2,k}(-b_{k}) \in \ESp_{2m}(R).$ Then 
                $V\overset{\ESp} \sim \begin{bmatrix}
                 1  & 0 & 0 & \cdots & 0\\
                 c & 1 & 0 &\cdots & 0\\
                \end{bmatrix} ~\mbox{for~some} ~c\in A . $ Now take 
                $\beta = se_{21}(-c)$, then we gets, $V\overset{\ESp} \sim \begin{bmatrix}
                 1  & 0 &  \cdots & 0\\
                 0 & 1 &\cdots & 0\\
                \end{bmatrix}.  $  Therefore, $V$ is completable to an elementary symplectic matrix.
                \par Now assume that $n>1$, since $R$ is a local ring, $V \in \Sp \Um_{2n,2m}(R) \subseteq \Um_{2n,2m}(R)$, there exists 
$\varepsilon \in \E_{2m}(R)$ such that 
$$ V\varepsilon = \begin{bmatrix}
                 1  & 0 &  \cdots& 0 & \cdots & 0\\
                 0 & 1 & \cdots &0&\cdots & 0\\
                    & &   \vdots\\
0& 0&\cdots &1& \cdots& 0
                \end{bmatrix}.$$
\par In view of Rao-Swan Lemma, there exists $\rho \in \E_{2m-1}(R)$ such that $\varepsilon (1\perp \rho) \in \ESp_{2m}(R),$
 therefore 
 $$ V\varepsilon (1\perp \rho) = \begin{bmatrix}
                 1  & 0 &  \cdots & 0\\
                 0 & b_{2} & \cdots & b_{2m}\\
\hspace{0.5in}& W&\hspace{0.5in}\\
                \end{bmatrix} \in \Sp \Um_{2n,2m}(R),$$ $~\mbox{for~some}~b_{i} \in R, ~2\leq i \leq 2m, W \in \Sp \Um_{2n-2,2m}(R).$
\par
Repeating the process done in $n=1$ case,  there exists $\varepsilon_{1} \in \ESp_{2m}(R)$ such that 
$$V\varepsilon_{1}
 = \begin{bmatrix}
                 1  & 0 &  \cdots & 0\\
                 0 & 1 & \cdots & 0\\
\hspace{0.5in}& V'&\hspace{0.5in} \\
                \end{bmatrix} ~\mbox{for~some}~V'\in \Sp \Um_{2(n-1),2m}(R).$$ Since $V\varepsilon_{1} \in \Sp \Um_{2n, 2m}(R),$ $ (V\varepsilon_{1})\psi_{m}(V\varepsilon_{1})^{t} = \psi_{n}.$ Therefore upon comparing the coefficients on the both side of the equation, one 
gets $ V' = (0, V'')$ for some $V''\in \Sp \Um_{2(n-1), 2(m-1)}(R).$
 \par By induction hypothesis $V'$ is completable to an elementary symplectic matrix, therefore $V$ is completable to an 
 elementary symplectic matrix.
 $~~~~~~~~~~~~~~~~~~~~~~~~~~~~~~~~~~~~~~~~~~~~~~~~~~~~~~~~~~~~~~~~~~~~~~~~~~~~~~~~~~~~~~~~~~~~~~~~~~~~~~~~~~~~~~\qedwhite$  
 
 \begin{prop}
\label{localgeneralisedhomotopysymplectic}
 Let $R$ be a local ring and $V \in \Sp \Um_{2n,2m}(R)$ for $m>n\geq 2$ or $m = n\geq 3.$
  Let $\delta \in \Sp_{2n}(R)$ be symplectic homotopic to identity. 
 Let $\delta(T)$ be a homotopy of $\delta.$ Then there exists,
 $ \sigma(T) \in \Sp_{2m}(R[T])$ with $\sigma(0) = \textit{Id}$ and $\sigma(T)^{-1}(\delta(T) \perp I_{2m-2n}) 
 \in \ESp_{2m}(R[T])$ 
 such that $$\delta(T)V = V\sigma(T).$$
\end{prop}
${\pf}$ In view of Lemma \ref{localsymplecticcompletable}, $V$ is completable to a matrix $W \in \Sp_{2m}(R)$.
Since $R$ is a local ring, 
$W\in \ESp_{2m}(R).$ By $(${\cite[Corollary 1.11]{kop}}$)$, $\ESp_{2m}(R[T]) \trianglelefteq \Sp_{2m}(R[T])$, for $m\geq 3,$ 
there exists 
$\varepsilon_{1}(T) \in \ESp_{2m}(R[T])$ such that 
$$(\delta(T)\perp I_{2m-2n})W(\delta(T)\perp I_{2m-2n})^{-1} = \varepsilon_{1}(T).$$
 Thus we have $(\delta(T)\perp I_{2m-2n})W= 
\varepsilon_{1}(T)W^{-1}W(\delta(T)\perp I_{2m-2n}).$ Again by normality of $\ESp_{2m}(R[T])$ in $\Sp_{2m}(R[T]$ for 
$m\geq 3,$ there exists $\varepsilon(T) \in \ESp_{2m}(R[T])$  such that 
 $$(\delta(T)\perp I_{2m-2n})W = W(\delta(T)\perp I_{2m-2n})\varepsilon(T).$$
Note that $\varepsilon(0) = \textit{Id}.$
Upon taking $\sigma(T) = (\delta(T)\perp I_{2m-2n})\varepsilon(T)$ and multiplying above equation by$ \begin{bmatrix}
                 I_{2n} & 0\\
                 0 & 0\\
                \end{bmatrix}$, we gets desired result.
$~~~~~~~~~~~~~~~~~~~~~~~~~~~~~~~~~~~~~~~~~~~~~~~~~~~~~~~~~~~~~~~~~~~~~~~~~~~~~~~~~~~~~~~~~~~~~~~~~~~~~~~~~~~~~~\qedwhite$

\begin{theo}
\label{generalisedhomotopysymplectic}
 Let $R$ be a commutative ring and $V\in \Sp \Um_{2n,2m}(R)$ with $m> n \geq 2 ~ \mbox{or}~ m = n\geq3.$ 
 Let $\delta \in \Sp_{2n}(R)$ be symplectic homotopic to identity. 
 Let $\delta(T)$ be a homotopy of $\delta.$
 Then $\exists ~\sigma(T) \in \Sp_{2m}(R[T],(T))$ such that
 $$\delta(T) V = V\sigma(T)~ \mbox{and}~ \sigma(T)^{-1}(\delta(T) \perp I_{2m-2n}) \in \ESp_{2m}(R[T],(T)).$$
 Moreover
 , if $\sigma(1) = \sigma$, then we have
 $\delta V = V\sigma$ and $\sigma^{-1}(\delta \perp I_{2m-2n}) \in \ESp_{2m}(R).$
\end{theo}
${\pf}$  Define,
\begin{align*}J = &\{s\in R~ | ~\delta(T)_{s}V_{s} = V_{s}\sigma(T)~ \mbox{for ~some}~
 \sigma(T) \in \Sp_{2m}(R_{s}[T], (T))~\\
&\mbox {with}~
\sigma(T)^{-1}(\delta(T)_{s}\perp I_{2m-2n}) \in \ESp_{2m}(R_{s}[T])\}.
\end{align*}
$\bf{Claim:}$ $J$ is an ideal.
\par
 For $s\in J, \lambda \in R$, clearly $\lambda s \in J$.
 So we need to prove that if $s_{1}, s_{2} \in J$ then
 $s_{1} + s_{2} \in J.$ Since $s_{1}, s_{2} \in J$,
we have $(s_{1} + s_{2})s_{1}$, $(s_{1} + s_{2})s_{2} \in J.$
We rename $R_{s_{1}+s_{2}}$ by $R$, now
it suffices to show that
\begin{align*}\delta(T)V &= V\sigma(T) ~\mbox{for}~\mbox{some}~ \sigma(T)
\in \Sp_{2m}(R[T], (T))\\
~&\mbox{with}~ 
\sigma(T)^{-1}(\delta(T)\perp I_{2m-2n})\in \ESp_{2m}(R[T], (T))~
\mbox{provided that}~ s_{1} + s_{2} = 1~\mbox{ and}
\end{align*}
\begin{equation}
\label{eq5}
\delta(T)_{s_{1}}V_{s_{1}} = V_{s_{1}}\sigma_{1}(T) ~\mbox{with}~
\sigma_{1}(T)^{-1}(\delta(T)_{s_{1}}\perp
I_{2m-2n})\in \ESp_{2m}(R_{s_{1}}[T], (T)),
\end{equation}
\begin{equation}
\label{eq6}
 \delta(T)_{s_{2}}V_{s_{2}} = V_{s_{2}}\sigma_{2}(T) ~\mbox{with}~ 
 \sigma_{2}(T)^{-1}(\delta(T)_{s_{2}}\perp
I_{2m-2n})\in \ESp_{2m}(R_{s_{2}}[T], (T)).
\end{equation}
Let  $$~~~~~~~~~~~~~~~~~~~~~~~~~~~~\sigma_{1}(T)(\delta(T)_{s_{1}} \perp
I_{2m-2n})^{-1} = \varepsilon_{1}(T) 
\in \ESp_{2m}(R_{s_{1}}[T], (T)),$$
       $$~~~~~~~~~~~~~~~~~~~~~~~~~~~~\sigma_{2}(T)(\delta(T)_{s_{2}} \perp 
       I_{2m-2n})^{-1} = \varepsilon_{2}(T) 
       \in \ESp_{2m}(R_{s_{2}}[T], (T)).$$

 Now, $\begin{bmatrix}
                 V_{s_{1}s_{2}}\\
                  0\\
                \end{bmatrix} \varepsilon_{1}(T)_{s_{2}}
                \varepsilon_{2}(T)_{s_{1}}^{-1}
                = \begin{bmatrix}
                 V_{s_{1}s_{2}}\\
                  0\\
                \end{bmatrix}.$ Let $\theta (T) =
                \varepsilon_{1}(T)_{s_{2}}\varepsilon_{2}(T)_{s_{1}}^{-1}.$ 
By Quillen's splitting property, $\mbox{for}~b\in (s_{2}^{N}),~ N>>0$, we have 
\begin{equation}
 \label{eq7}
 \theta(T) = \theta(bT)\{ \theta(bT)^{-1}\theta(T)\}
\end{equation} $~\mbox{with}
 ~\theta(bT)\in \ESp_{2m}(R_{s_{1}}[T]), ~ \mbox{and}~ 
 \theta(bT)^{-1}\theta(T)\in  \ESp_{2m}(R_{s_{2}}[T]).$

Since $\begin{bmatrix}
                 V_{s_{1}s_{2}}\\
                  0\\
                \end{bmatrix} \theta(T) = \begin{bmatrix}
                 V_{s_{1}s_{2}}\\
                  0\\
                \end{bmatrix}$. We have, $\begin{bmatrix}
                 V_{s_{1}s_{2}}\\
                  0\\
                \end{bmatrix} \theta(bT) = \begin{bmatrix}
                 V_{s_{1}s_{2}}\\
                  0\\
                \end{bmatrix}$ and \\ $\begin{bmatrix}
                 V_{s_{1}s_{2}}\\
                  0\\
                \end{bmatrix} \theta(bT)^{-1}\theta(T)
                 = \begin{bmatrix}
                 V_{s_{1}s_{2}}\\
                  0\\
                \end{bmatrix}.$
                \par
                Define, $\theta(bT) = \eta_{1}(T)^{-1}, ~ \theta(bT)^{-1}\theta(T) = \eta_{2}(T)$. Thus we have 
                $V_{i}\eta_{i}(T) = V_{i}$ with $\eta_{i}(0) = \textit{Id}$ and $\eta_{i}(T) \in \ESp_{2m}(R_{s_{i}}[T])$ 
                for $i =1,2.$
                \par
                In view of equation \ref{eq7}, we gets 
                \begin{equation}
 \label{eq8}
 (\eta_{1}(T)\varepsilon_{1}(T))_{s_{2}} = (\eta_{2}(T)\varepsilon_{2}(T))_{s_{1}}.
\end{equation}
Now, by equations \ref{eq5} and \ref{eq6}, 
$$\delta(T)_{s_{1}}V_{s_{1}} = V_{s_{1}}\eta_{1}(T)\sigma_{1}(T)$$
$$\delta(T)_{s_{2}}V_{s_{2}} = V_{s_{2}}\eta_{2}(T)\sigma_{2}(T).$$
In view of equation \ref{eq8}, we have $(\eta_{1}(T)\sigma_{1}(T))_{s_{2}} = (\eta_{2}(T)\sigma_{2}(T))_{s_{1}}$. Since 
$s_{1} + s_{2} = 1$, 
$\exists ~\sigma(T) \in \Sp_{2m}(R[T])$ such that $\sigma(T)_{s_{1}} = \eta_{1}(T)\sigma_{1}(T)$ and 
$\sigma(T)_{s_{2}} = \eta_{2}(T)\sigma_{2}(T)$ with $\sigma(T)_{s_{i}}^{-1}(\delta(T)_{s_{i}}\perp I_{2m-2n})
\in \ESp_{2m}(R_{s_{i}}
[T])$ for $i = 1, 2.$ Since $s_{1}$ and $s_{2}$ are comaximal, 
by theorem \ref{2.17},  
\begin{align*}\delta(T)V &= V\sigma(T) ~\mbox{for}~\mbox{some}~ 
 \sigma(T) \in \Sp_{2m}(R[T], (T)) ~\\
&\mbox{with}~ 
\sigma(T)^{-1}(\delta(T)\perp I_{2m-2n})\in \ESp_{2m}(R[T], (T)).
\end{align*}
This proves that $J$ is an ideal.
\par
In view of lemma \ref{localgeneralisedhomotopysymplectic}, for every maximal ideal $\mathfrak{m}$ of $R$, we have 
$$\delta(T)_{\mathfrak{m}}V_{\mathfrak{m}} = V_{\mathfrak{m}}\sigma^{'}(T) ~\mbox{with}~ 
\sigma^{'}(T)^{-1}(\delta(T)_{\mathfrak{m}}\perp
I_{2m-2n})\in \ESp_{2m}(R_{\mathfrak{m}}[T], (T)).$$
Thus there exists $s \in R\setminus \mathfrak{m}$, such that
$$\delta(T)_{s}V_{s} = V_{s}\sigma^{'}(T) ~\mbox{with}~ 
\sigma^{'}(T)^{-1}(\delta(T)_{s}\perp
I_{2m-2n})\in \ESp_{2m}(R_{s}[T], (T)).$$
Therefore $J \nsubseteq \mathfrak{m}$, for any maximal ideal $\mathfrak{m}$ of $R$ i.e. $1 \in J.$
Thus $\exists ~\sigma(T) \in \Sp_{2m}(R[T], (T))$ such that 
$$\delta(T)V = V\sigma(T) ~\mbox{with}~ 
\sigma(T)^{-1}(\delta(T)\perp I_{2m-2n})\in \ESp_{2m}(R[T], (T)).$$
Now put $T=1$, and take $\sigma(1) = \sigma$ to get the desired result.
$~~~~~~~~~~~~~~~~~~~~~~~~~~~~~~~~~~~~~~~~~~~~~~~~~~~~~~~~~~~~~~~~~~~~~~~~~~~~~~~~~~~~~~~~~~~~~~~~~~~~~~~~~~~~~~\qedwhite$

 \begin{cor}
 $(${\cite[Theorem 2.19]{rs}}$)$
 Let $m\geq 2$ and $\alpha, \beta \in \Sp_{2m}(R).$  Let either $\alpha$ or $\beta$ be symplectic homotopic to identity.
 Then $\alpha \beta = \beta \alpha \varepsilon,~
  ~\mbox{for~ some} ~\varepsilon \in \ESp_{2m}(R).$
\end{cor}
${\pf}$ Let us assume that $\alpha$ is homotopic to identity, so there exists
 $\delta(T) \in \Sp_{2m}(R[T])$ such that $\delta(0) = \textit{Id}~\mbox{and}~ \delta(1) = \alpha.$ By theorem
 \ref{generalisedhomotopysymplectic}, 
 there exists $\varepsilon(T) \in \ESp_{2m}(R[T])$ with $\varepsilon(0) = \textit{Id}$ such that 
 $$\delta(T)\beta = \beta\delta(T)\varepsilon(T).$$
 \par Put $T =1$ to get the desired result.
 $~~~~~~~~~~~~~~~~~~~~~~~~~~~~~~~~~~~~~~~~~~~~~~~~~~~~~~~~~~~~~~~~~~~~~~~~~~~~~~~~~~~~~~~~~~~~~~~~~~~~~~~~~~~~~~\qedwhite$
                
\begin{cor}
 Let $\delta \in \Sp_{2n}(R)$ and $V\in \Sp \Um_{2n,2m}(R)$. Then $\delta V = V\sigma$ for some $\sigma \in \Sp_{2m}(R)$ such that
  $(\delta^{-1} \perp \sigma) \in \ESp_{2(n+m)}(R).$
\end{cor}
${\pf}$ By $(${\cite[Lemma 1.1]{vas1}}$)$, $(\delta \perp \delta^{-1}) \in \ESp_{4n}(R).$ Since every elementary 
symplectic matrix is homotopic 
to identity, 
thus by theorem \ref{generalisedhomotopysymplectic}, 
$$(\delta \perp \delta^{-1})(V\perp I_{2n}) = (V\perp I_{2n})\sigma^{'}, ~\mbox{with}~\sigma^{'} \in \ESp_{2(n+m)}(R).$$
Write $\sigma^{'} = \begin{bmatrix}
                 \alpha & \beta\\
                 \gamma & \zeta\\
                \end{bmatrix}$ where $\alpha \in M_{2m\times 2m}(R),~\beta\in M_{2m\times 2n}(R), ~ 
                \gamma \in M_{2n\times 2m}(R),
                ~ \zeta \in M_{2n\times 2n}(R)$. Thus we have, 
                $$(\delta V \perp \delta^{-1}) = \begin{bmatrix}
                 V\alpha & V\beta\\
                 \gamma& \zeta\\
                \end{bmatrix}.$$
 \par
 Upon compairing both sides we gets $\gamma= 0$ and $\zeta = \delta^{-1}.$ Therefore\\
 $$\begin{bmatrix}
                 \alpha & \beta\\
                 0 & \delta^{-1}\\
                \end{bmatrix} \in \ESp_{2(n+m)}(R).$$
\par
Now, take $\alpha = \sigma$, so we have $(\delta^{-1} \perp \sigma) \in \ESp_{2(n+m)}(R)$ and $\delta V = V\sigma.$
$~~~~~~~~~~~~~~~~~~~~~~~~~~~~~~~~~~~~~~~~~~~~~~~~~~~~~~~~~~~~~~~~~~~~~~~~~~~~~~~~~~~~~~~~~~~~~~~~~~~~~~~~~~~~~~\qedwhite$

\section{\bf{Generalised Homotopy and Commutativity Principle for Orthogonal Groups}}                
\par
  Throughout this section we will assume that $1/2 \in R$, where $R$ is a commutative ring with $1\neq 0.$  In this section, we will deal with orthogonal matrices of size at least 6.
                
\begin{notn}

Let  
                $\phi_{1} = \begin{bmatrix}
                 0 & 1\\
                  1 & 0\\
                \end{bmatrix},~~ \phi_{n} =\phi_{n-1} \perp \phi_{1};$ for $n > 1$.
\end{notn}

\begin{notn}
 Let $\sigma$ be the permutation of the natural numbers given by $\sigma(2i) = 2i-1$ and $\sigma(2i-1) = 2i$.
\end{notn}

   \begin{defi}{\bf{Orthogonal group}}~{$\rm{O}_{2m}(R) :$} The subgroup of $\GL_{2m}(R)$ consisting of all $2m \times 2m $ matrices 
$\{\alpha \in \GL_{2m}(R)~\mid \alpha^{t}\phi_{m}\alpha = \phi_{m}\}$.
\end{defi}
\begin{defi}{\bf{Elementary orthogonal group}}~{$\EO_{2m}(R) :$} We define for $ 1\leq i \neq j\leq 2m,~z\in R,$\\
$$oe_{ij}(z) = I_{2m} + zE_{ij} - zE_{\sigma(j)\sigma(i)}, ~~~\textit{if} ~ i\neq \sigma(j).$$\\
It is easy to verify that all these matrices belong to $\rm{O}_{2m}(R)$. We call them the elementary orthogonal matrices over
$R$. The subgroup generated by them is called the elementary orthogonal group and is denoted by $\EO_{2m}(R)$.
\end{defi}

   \begin{defi}
 $\mbox{OUm}_{2n, 2m}(R) = \{V \in \Um_{2n,2m}(R) | V\phi_{m}V^{t} = \phi_{n}\}.$
\end{defi}

\par Let $P$ be a finitely generated projective $R$-module. The module $P \oplus P^{\ast}$ carries a natural quadratic form $q$ defined by $q(x+f) = f(x)$ for $x\in P$ and $f\in P^{\ast}.$ The associated bilinear form is given by $B_{q}(x_{1}+f_{1}, x_{2}+f_{2}) = 
f_{1}(x_{2}) + f_{2}(x_{1}), x_{1}, x_{2}\in P, f_{1}, f_{2} \in P^{\ast}.$  It is easy to see that $q$ is non-singular. The quadratic space $(P\oplus P^{\ast}, q)$ will be called hyperspace of $P.$ The hyperbolic space of a free $R$-module of rank 1 is called a hyperplane. 
\begin{defi} An orthogonal pair of elements $(w_{1}, w_{2})$ is said to be a hyperbolic pair if $q(w_{1}) = 1, q(w_{2}) = -1.$
\end{defi}

\begin{rem} Hyperbolic plane is generated as an $R$-module by a hyperbolic pair.
\end{rem}

\begin{lem}
\label{orthocomplete}
 Let $R$ be a local ring with $2R = R$ and $V\in \rm{OUm}_{2n,2m}(R), m \geq  n+2, n \geq 1.$ Then $V$ is completable to an elementary orthogonal matrix.
\end{lem}

${\pf}$ We will prove it by induction on $n, m.$ Let us assume that $n = 1.$ In view of $(${\cite[Theorem 7.1 (ii)]{roy}}$)$ and  $(${\cite[Lemma 2.7]{ambily}}$),$  there exists $\varepsilon \in \EO_{2m}(R)$ such that 
$$V\varepsilon = \begin{bmatrix}
                 1  & 0 &  \cdots & 0\\
                 0 & 1 & \cdots & 0\\
                \end{bmatrix}.$$
Thus $V$ is completable to an elementary orthogonal matrix.
 \par  Now assume that $n>1$. We have $m > 3.$ In view of $(${\cite[Theorem 7.1 (ii)]{roy}}$),$ there exists $\varepsilon_{1} \in \EO_{2m}(R)$ such that 
$$V\varepsilon_{1}
 = \begin{bmatrix}
                 1  & 0 &  \cdots & 0\\
                 0 & 1 & \cdots & 0\\
\hspace{0.5in}& V'&\hspace{0.5in} \\
                \end{bmatrix} ~\mbox{for~some}~V'\in \rm{OUm}_{2(n-1),2m}(R).$$ Since $V\varepsilon_{1} \in \mbox{OUm}_{2n, 2m}(R),$ $ (V\varepsilon_{1})\phi_{m}(V\varepsilon_{1})^{t} = \phi_{n}.$ Therefore upon comparing the coefficients on the both side of the equation, one 
gets $ V' = (0, V'')$ for some $V''\in \rm{OUm}_{2(n-1), 2(m-1)}(R).$ Now, we get the desired result by induction hypothesis.
$~~~~~~~~~~~~~~~~~~~~~~~~~~~~~~~~~~~~~~~~~~~~~~~~~~~~~~~~~~~~~~~~~~~~~~~~~~~~~~~~~~~~~~~~~~~~~~~~~~~~~~~~~~~~~~\qedwhite$

Following the steps of the proof of Proposition \ref{localgeneralisedhomotopysymplectic}, one gets the following result :

 \begin{prop}
Let $R$ be a local ring and $V \in \rm{OUm}_{2n,2m}(R)$ for $m\geq n+2, n \geq 2.$
  Let $\delta \in \rm{SO}_{2n}(R)$ be orthogonal homotopic to identity. 
 Let $\delta(T)$ be a homotopy of $\delta.$ Then there exists,
 $ \sigma(T) \in \mbox{SO}_{2m}(R[T])$ with $\sigma(0) = \textit{Id}$ and $\sigma(T)^{-1}(\delta(T) \perp I_{2m-2n}) 
 \in \EO_{2m}(R[T])$ 
 such that $$\delta(T)V = V\sigma(T).$$
\end{prop}

 By making appropriate modifications in the proof of theorem \ref{generalisedhomotopylinear} and theorem \ref{generalisedhomotopysymplectic}, one can prove the following result : 
\begin{theo}
\label{generalisedhomotopyorthogonal}
 Let $R$ be a commutative ring and $V\in \rm{OUm}_{2n,2m}(R)$ with $m\geq n+2, n \geq 2 .$  Let $\delta \in \rm{SO}_{2n}(R)$ be orthogonal homotopic to identity. 
 Let $\delta(T)$ be a homotopy of $\delta.$
 Then $\exists ~\sigma(T) \in \mbox{SO}_{2m}(R[T],(T))$ such that
 $$\delta(T) V = V\sigma(T)~ \mbox{and}~ \sigma(T)^{-1}(\delta(T) \perp I_{2m-2n}) \in \EO_{2m}(R[T],(T)).$$
 Moreover
 , if $\sigma(1) = \sigma$, then we have
 $\delta V = V\sigma$ and $\sigma^{-1}(\delta \perp I_{2m-2n}) \in \EO_{2m}(R).$
\end{theo}

\par Due to the size restrictions in lemma \ref{orthocomplete} one is not able to deduce whether a similar homotopy and commutativity principle holds in the orthogonal case. We began this study in \cite{rs}. We add a few more observations on this below.

\begin{lem}
 \label{vasraolike}
 $($L.N. Vaserstein$)$ $(${\cite[Theorem 3.5]{vas1}}$)$ Let $m\geq 3$ and $R$
 be a local ring, $\frac{1}{2}\in R$. Then ${\rm{O}_{2m}(R)}/{\EO_{2m}(R)}
 = {\rm{O}_{2}(R)}/{\EO_{2}(R)} = O_{2}(R).$
 \end{lem}
 \begin{obs}
  \label{orthoobservation} Every element $\alpha \in \rm{O}_{2}(R)$ is either of the type $ \begin{bmatrix}
                 u  & 0 \\
                 0 &  u^{-1}\\
                \end{bmatrix}$ or of the type $ \begin{bmatrix}
                 0  & u \\
                 u^{-1} &  0\\
                \end{bmatrix}$ for some $u\in R^{\ast}.$
 \end{obs}

\begin{theo}
\label{ortho1stably}
Let $R$ be a local ring, $m\geq 3$ and $1/2 \in R$. Then we have, $$([\rm{O}_{2m}(R[X]), \rm{O}_{2m}(R[X])]\perp I_{2}) \subseteq \EO_{2m+2}(R[X]).$$
\end{theo}
${\pf}$ Let $\alpha(X), \beta(X) \in \rm{O}_{2m}(R[X])$, we need to prove that $([\alpha(X), \beta(X)] \perp I_{2}) \in
\EO_{2m+2}(R[X]).$ Define,
$$\gamma(X, T) = [\alpha(XT)\perp I_{2}, \beta(X) \perp I_{2}]$$
\par For every maximal ideal $\mathfrak{m}$ of $R[X]$, we have 
$\gamma(X, T)_{\mathfrak{m}} = [(\alpha(XT)\perp I_{2})_{\mathfrak{m}}, (\beta(X) \perp I_{2})_{\mathfrak{m}}].$ In view of 
lemma \ref{vasraolike}, $(\beta(X) \perp I_{2})_{\mathfrak{m}} = (I_{2m}\perp \delta(X))\varepsilon(X)$ 
for some $\delta(X) \in \rm{O}_{2}(R[X]_{\mathfrak{m}})$ and $\varepsilon(X) \in \EO_{2m+2}(R[X])_{\mathfrak{m}}.$ 
By observation \ref{orthoobservation}, either $\delta(X) = \begin{bmatrix}
                 u & 0\\
                 0 & u^{-1}\\
                \end{bmatrix}$\\ or $\delta(X) = \begin{bmatrix}
                 0 & u\\
                 u^{-1} & 0\\
                \end{bmatrix},~\mbox{for~some}~u\in R[X]_{\mathfrak{m}}^{\ast}.$ Therefore 
$\gamma(X, T)_{\mathfrak{m}} \in \EO_{2m+2}(R[X]_{\mathfrak{m}}[T]).$
\par Now, $\gamma(X,0) = [\alpha(0)\perp I_{2}, \beta(X) \perp I_{2}].$ Since $R$ is a local ring, by lemma
 \ref{vasraolike}, $\alpha(0) \perp I_{2} = 
(I_{2m} \perp \delta)\varepsilon_{1}$ for $\delta \in \rm{O}_{2}(R)$ and $\varepsilon_{1} \in \EO_{2m+2}(R).$ By 
observation \ref{orthoobservation}, either $\delta = \begin{bmatrix}
                 a & 0\\
                 0 & a^{-1}\\
                \end{bmatrix}~or~ \begin{bmatrix}
                 0 & a\\
                 a^{-1} & 0\\
                \end{bmatrix},~\mbox{for~some}~a\in R^{\ast}.$ Therefore, 
$\gamma(X,0) \in \EO_{2m+2}(R[X]).$ Now by local-global principle for othogonal groups $(${\cite[Theorem 4.2] {sus}}$)$, 
we have 
$$\gamma(X, 1) = [\alpha(X)\perp I_{2}, \beta(X) \perp I_{2}] = ([\alpha(X), \beta(X)] \perp I_{2}) \in \EO_{2m+2}(R[X]).$$
$~~~~~~~~~~~~~~~~~~~~~~~~~~~~~~~~~~~~~~~~~~~~~~~~~~~~~~~~~~~~~~~~~~~~~~~~~~~~~~~~~~~~~~~~~~~~~~~~~~~~~~~~~~~~~~\qedwhite$

\begin{notn}
 We will denote set of all special orthogonal matrices which are special orthogonally homotopic to identity by $\rm{HSO}_{2m}(R).$
\end{notn}

\begin{theo}
 Let $m\geq 2$ and $R$ be a commutative ring,  $\frac{1}{2}\in R$. Then, $$[\rm{HSO}_{2m}(R)\perp I_{2}, \rm{O}_{2m}(R)\perp I_{2}] \subseteq \EO_{2m+2}(R).$$
\end{theo}

${\pf}$ Let $\alpha \in \rm{HSO}_{2m}(R), \beta \in \rm{O}_{2m}(R)$, we need to prove that $[\alpha \perp I_{2}, \beta\perp I_{2}] 
 \in \EO_{2m+2}(R).$ Let $\alpha(T)$ be a homotopy of $\alpha$ and define, 
 $$\gamma(T) = [\alpha(T) \perp I_{2}, \beta \perp I_{2}].$$
 \par Clearly, $\gamma(0) = Id.$ For every maximal ideal $\mathfrak{m}$ of $R$, we have 
$$\gamma(T)_{\mathfrak{m}} = [(\alpha(T) \perp I_{2})_{\mathfrak{m}}, (\beta \perp I_{2})_{\mathfrak{m}}].$$ 
In view of lemma
 \ref{vasraolike}, $(\beta \perp I_{2})_{\mathfrak{m}} = 
(I_{2m} \perp \delta)\varepsilon$ for some $\delta \in \rm{O}_{2}(R_{\mathfrak{m}})$ and $\varepsilon \in 
\EO_{2m+2}(R_{\mathfrak{m}}).$
By observation \ref{orthoobservation}, either $\delta = \begin{bmatrix}
                 a & 0\\
                 0 & a^{-1}\\
                \end{bmatrix}~or~ \begin{bmatrix}
                 0 & a\\
                 a^{-1} & 0\\
                \end{bmatrix},~\mbox{for~some}~a\in R_{\mathfrak{m}}^{\ast}.$ Therefore,
                 $\gamma(T)_{\mathfrak{m}} \in \EO_{2m}(R_{\mathfrak{m}}[T]).$ In view of local-global principle for 
                orthogonal groups  $(${\cite[Theorem 4.2] {sus}}$)$, we have
                $\gamma(T) \in \EO_{2m+2}(R[T]).$ Therefore, $\gamma(1) = [\alpha \perp I_{2}, \beta\perp I_{2}]\in 
                \EO_{2m+2}(R).$
$~~~~~~~~~~~~~~~~~~~~~~~~~~~~~~~~~~~~~~~~~~~~~~~~~~~~~~~~~~~~~~~~~~~~~~~~~~~~~~~~~~~~~~~~~~~~~~~~~~~~~~~~~~~~~~\qedwhite$

\medskip
\noindent
{\bf Acknowledgement:} The first author thanks Professor L.N. Vaserstein for some valuable insights. The second author is thankful to Inspire Faculty Fellowship (DST/INSPIRE/04/2021/002849) and Startup Research Grant (SRG/2022/000056) for their support. The second author also thanks IIT Mandi for their Seed Grant.

\Addresses

\end{document}